\documentclass[11pt, a4paper]{article}
\usepackage{amsmath, amsfonts, amssymb, amscd, amsthm, mathdots, mathrsfs, dsfont}
\usepackage{titlesec}
\usepackage{fancyhdr}
\usepackage{fullpage}
\usepackage[colorlinks, urlcolor=black]{hyperref}
\hypersetup{citecolor=blue, linkcolor=blue}
\usepackage{tikz}
\usepackage{mathabx}

\newcommand{\N}{\mathbb{N}}
\newcommand{\Q}{\mathbb{Q}}
\newcommand{\R}{\mathbb{R}}
\newcommand{\Z}{\mathbb{Z}}

\DeclareMathOperator{\Disc}{Disc}
\DeclareMathOperator{\Gal}{Gal}

\titleformat{\section}[hang]
{\normalfont\filright\large}{\thesection . }{0pt}
{\upshape\bfseries}

\titleformat{\subsection}[hang]
{\itshape}{\thesubsection \ - }{0pt}
{}

\newtheorem{theorem}{Theorem}[section]
\newtheorem{proposition}[theorem]{Proposition}
\newtheorem{conjecture}[theorem]{Conjecture}
\newtheorem{corollary}[theorem]{Corollary}
\newtheorem{lemma}[theorem]{Lemma}

\theoremstyle{definition}

\newtheorem{question}{Question}

\theoremstyle{remark}
\newtheorem{remark}{Remark}

\title{On the inverse Galois problem for domination polynomials}
\author{\small JOSÉ LUIS COSME-ÁLVAREZ \& ADRIÁN ZENTENO}
\date{\today}
\begin{document}

\maketitle

\begin{abstract}

The classical inverse Galois problem asks whether every finite group occurs as the Galois group of a polynomial over $\mathbb{Q}$. A more restrictive problem arise when the defining polynomial is required to come from a prescribed combinatorial construction. 

In this paper we investigate this problem for domination polynomials of graphs. We show that complete graphs realize the groups $(\mathbb{Z}/n\mathbb{Z})^\times$, while stars realize $S_n$ for $n\not\equiv5\pmod6$ and $\mathbb{Z}/2\mathbb{Z}\times S_{n-2}$ for $n\equiv5\pmod6$. Oboudi's factorization of sunlike graphs then propagates these realizations to infinitely many trees. We place this construction in a more general framework based on irrelevant edges between domination-covered vertices, which provides connected realizations of products and yields, among other examples, $(\mathbb{Z}/2\mathbb{Z})^4$ and $D_5\times\mathbb{Z}/2\mathbb{Z}$.
Then, we relate the inverse Galois problem to Oboudi's conjecture on domination polynomials with real roots: using the fact that a polynomial with Galois group of odd order has only real roots, we prove that Oboudi's conjecture would rule out the realization of any nontrivial finite group of odd order as the Galois group of a domination polynomial. In particular, realizing $\mathbb{Z}/3 \mathbb{Z}$, or indeed any nontrivial finite group of odd order, would disprove the conjecture.
Finally, for comparison, we prove that the elementary abelian $3$-groups $(\mathbb{Z}/3\mathbb{Z})^m$ occur as Galois groups of chromatic polynomials for every $m\geq0$. More generally,  we prove that every chromatically realizable group $\Gamma$ gives connected realizations of $\Gamma\times (\mathbb{Z}/2\mathbb{Z})^n\times (\mathbb{Z}/3\mathbb{Z})^m$ for all $n,m\geq0$.

\end{abstract}

%\tableofcontents

\section{Introduction}\label{sec:intro}

One of the basic constructions in Galois theory associates a finite group to a polynomial. Let $f(x)\in \Q[x]$. We say that $f(x)$ splits in a field extension $L/\Q$ if every root of $f(x)$ is contained in $L$. The splitting field $\Q_{f(x)}$ of $f(x)$ is the smallest extension of $\Q$ in which $f(x)$ splits, and its Galois group is  
\[
\Gal(f(x)) := \Gal(\Q_{f(x)}/\Q),
\]
the group of automorphisms of $\Q_{f(x)}$. If $f(x)$ has $n$ distinct roots, $\Gal(f(x))$ can be identified with a subgroup of the symmetric group $S_n$, because every automorphism of $\Q_{f(x)}$ permutes the roots of $f(x)$ and is determined by its action on them.
This construction naturally leads to the following question: given a finite group $\Gamma$, does there exist a polynomial over $\Q$ whose Galois group is isomorphic to $\Gamma$? This question is known as the inverse Galois problem. It is conjectured to have an affirmative answer, although it remains open in general. Among the earliest positive results, Hilbert \cite{Hilbert1892} showed that, for every $n$, the symmetric and alternating groups $S_n$ and $A_n$ occur as Galois groups over $\Q$. For a comprehensive account of these and related results, see \cite{MalleMatzat2018}.

Inspired by this classical problem, one may impose the additional requirement that the polynomial arise from a prescribed combinatorial construction. If $F(G,x)$ is a graph polynomial, the corresponding inverse Galois question is: 
\begin{quote}
\emph{Does every finite group occur as the Galois group of $F(G,x)$ for some finite simple graph $G$?}
\end{quote}
The restriction to graph polynomials changes the nature of the problem: the coefficients encode combinatorial information about the graph and are subject to structural constraints absent from arbitrary polynomials over $\Q$.
 
The chromatic polynomial $P(G,x)$ of a graph $G$ provides a natural first example of such a graph polynomial. The first systematic study of its Galois groups was carried out by Morgan \cite{Morgan2010, Morgan2012}. Her computations showed a predominance of symmetric groups among the irreducible factors of chromatic polynomials of small graphs, but she also constructed several infinite families with non-symmetric Galois groups, including infinite families whose chromatic polynomials have cyclic Galois groups of order 2 and 3. Further developments include the theta-graph analysis of Delbourgo–Morgan \cite{DelbourgoMorgan2014}, the biclique constructions of Gilmore \cite{DelbourgoMorgan2019, Gilmore2015}, and the broader algebraic study of Cameron–Morgan \cite{CameronMorgan2017}. Nevertheless, the chromatic inverse Galois problem remains largely open \cite[Problem 3]{Morgan2012}.
 In Section \ref{sec:3}, we revisit this problem, recalling some known realizations and proving new realization results involving elementary abelian 2- and 3-groups.  More precisely, we prove that, if a finite group $\Gamma$ is realizable as the Galois group of a chromatic polynomial, then the parameters in Morgan's quadratic and cyclic cubic families can be chosen so that the corresponding splitting fields have pairwise disjoint ramification and avoid the ramification of a fixed realization of $\Gamma$. Thus, suitable $K_1$-gluings yield a connected graph $G$ such that 
   \[
\Gal(P(G,x)) \cong  \Gamma\times(\Z/2\Z)^n\times(\Z/3\Z)^m
  \]
  for arbitrary $n,m\ge0$. In particular, combining this result with the Galois groups of equal-length theta graphs of Delbourgo--Morgan \cite[Theorem 1.3.i]{DelbourgoMorgan2014}, we obtain $(\Z/2\Z)^n\times(\Z/3\Z)^m\times S_{3k}$ as the Galois group of a chromatic polynomial for all $n,m\ge0$ and $k\ge1$. 
  
The chromatic case provides a useful point of comparison for our main object of study, the domination polynomial $D(G,x)$. In the form now standard in the literature, $D(G,x)$ was introduced by Alikhani and Peng in \cite{AlikhaniPengIntro2014}, although polynomials enumerating dominating sets had already been considered by Arocha and Llano in \cite{ArochaLlano2000}. Since then, domination polynomials have been studied from analytic and combinatorial viewpoints \cite{AlikhaniPeng2014, BeatonBrown2021, BrownTufts2014, Oboudi2016}. However, in contrast with the chromatic case, their Galois groups do not appear to have been studied previously.  
    
In Section \ref{sec:4}, we begin by exhibiting two natural infinite families of realizations. First, we prove that the domination polynomial $D(K_n,x)$ of a complete graph $K_n$ has cyclotomic splitting field $\Q(\zeta_n)$ and Galois group $(\Z/n\Z)^\times$, thereby realizing an infinite family of abelian groups. Stars $K_{1,n}$ exhibit a very different phenomenon. Combining a result of Osada \cite{Osada1987} with the analysis of $X^n+X+1$ by Delbourgo--Morgan \cite{DelbourgoMorgan2014}, we obtain
\[
 \Gal(D(K_{1,n},x))\cong
 \begin{cases}
   S_n,&n\not\equiv5\pmod 6,\\
   \Z/2\Z \times S_{n-2},&n\equiv5\pmod6.
 \end{cases}
\]
Thus, stars realize symmetric groups $S_n$ for five residue classes modulo 6, whereas in the chromatic setting, the equal-length theta-graph construction of Delbourgo--Morgan realizes $S_{3k}$.
We also show that our star construction extends naturally to sunlike graphs. Let $G$ be a graph with vertices $v_1,\ldots,v_r$. If we attach $m_i\ge1$ new leaves to $v_i$, we obtain the sunlike graph $G(v_1^{m_1},\ldots,v_r^{m_r})$, and a result of Oboudi \cite{Oboudi2016} gives
\begin{equation}\label{facto}
 D(G(v_1^{m_1},\ldots,v_r^{m_r}),x) =\prod_{i=1}^rD(K_{1,m_i},x).
\end{equation}
Consequently, attaching the same number $m$ of leaves to every vertex of a tree produces another tree whose domination polynomial has the same splitting field as $D(K_{1,m},x)$. Hence every group realized by a star admits infinitely many pairwise nonisomorphic tree realizations, including infinitely many caterpillars obtained by taking paths. 

The chromatic polynomial has a particularly useful product behavior under clique gluings, which we in fact use in Section \ref{sec:3} to prove our main result for chromatic polynomials. For domination polynomials, there is no direct analogue for arbitrary vertex identifications. However, we place the Oboudi's factorization \eqref{facto} in a more general framework based on a theory of irrelevant edges, developed by Kotek--Preen--Simon--Tittmann--Trinks \cite{KotekEtAl2012}, which provides a remarkably close substitute. When connected graphs possess suitable domination-covered vertices, they can be joined by irrelevant edges without changing the product of their domination polynomials. Whenever the corresponding splitting fields are linearly disjoint, this construction realizes the direct product of their Galois groups while preserving connectedness. As an application, we construct connected graphs with domination Galois groups $(\Z/2\Z)^4$ and $D_5\times\Z/2\Z$, which cannot arise from a complete or a star graph.

In contrast with the chromatic case, where Cameron and Morgan conjecture that every finite group is chromatically realizable \cite[Problem 3]{Morgan2012}, we show in Section \ref{sec:5} that a conjecture of Oboudi \cite[Conjecture 1]{Oboudi2016} concerning the characterization of graphs whose domination polynomials have only real roots would imply that the inverse Galois problem for domination polynomials has a negative answer. In fact, we prove that this conjecture implies that no nontrivial finite group of odd order is realizable as a Galois group of a domination polynomial. Then, the realization of $\Z/3\Z$, or more generally of any nontrivial finite group of odd order, would disprove Oboudi's conjecture. This gives a sharp contrast with the chromatic case, where we prove that elementary abelian \(3\)-groups of arbitrary rank are realizable.

Finally, in Section \ref{sec:6} we compare the chromatic and domination inverse Galois problems for groups of small order and for transitive groups of small degree. Using the known chromatic realizations and our results for domination polynomials, we show that the chromatic inverse Galois problem is affirmative for every group of order less than $6$, while the domination inverse Galois problem is affirmative for every even-order group of order less than $6$. In both settings, the smallest even-order groups for which realizability remains open are $\Z/8\Z$ and $Q_8$. On the other hand, we have that the work of Cameron--Morgan \cite{CameronMorgan2017}, together with that of Gilmore \cite{DelbourgoMorgan2019, Gilmore2015}, shows that the known chromatic constructions realize every transitive group of degree at most $5$, while their degree-$6$ computations realize $5$ of the $16$ transitive groups. In the domination setting, we realize every transitive group of degree at most $4$ except $\Z/3\Z \cong A_3$ and $A_4$, and we give explicit quintic realizations of $D_5$ and $S_5$.

%%%%%%%%%%%%%%%%%%%%%%%%%%%%%%%%%%%%%%%%%%%%%%%%%%%%%%%%%%%%%%%%%%%%
%%%%%%%%%%%%%%%%%%%%%%%%%%%%%%%%%%%%%%%%%%%%%%%%%%%%%%%%%%%%%%%%%%%%

\section{Algebraic preliminaries}\label{sec:2}

In this section, we recall some basic notions and results from Galois theory and number theory that will be used throughout the paper. We refer the reader to \cite{StewartGalois, StewartNumberTheory} for further background.

Let $f(x)\in\Q[x]$, and let $\alpha_1,\dots,\alpha_n$ be its different roots in a fixed algebraic closure $\overline{\Q}$. The splitting field of $f(x)$ is
\[
\Q_{f(x)}:=\Q(\alpha_1,\dots,\alpha_n),
\]
the smallest extension of $\Q$ in $\overline{\Q}$ containing all the roots of $f(x)$, and its Galois group is 
\[
\Gal(f(x)) := \Gal(\Q_{f(x)}/\Q),
\] 
the group of field automorphisms of $\Q_{f(x)}$, which in particular fix $\Q$. Since every automorphism of $\Q_{f(x)}$ is determined by its action on the roots, $\Gal(f(x))$ may be viewed as a subgroup of the symmetric group $S_n$. If $f(x)$ is irreducible over $\Q$, then this action is transitive. For a reducible polynomial, the irreducible factors determine the orbits of the action. Factors of $f(x)$ that split completely over $\Q$ do not enlarge its splitting field. Thus, rational linear factors may be discarded when computing $\Gal(f(x))$.

If $\displaystyle f(x)=a_n\prod_{i=1}^n(x-\alpha_i)$, its discriminant is defined as
\[
  \Disc(f(x)) : =  a_n^{2n-2}  \prod_{1\leq i<j\leq n}  (\alpha_i-\alpha_j)^2.
\]
The discriminant of $f(x)$ provides useful information about the Galois group. For example, if $f(x)$ is separable of degree $n$, then $\Gal(f(x))$, viewed as a subgroup of $S_n$, is contained in $A_n$ if and only if $\Disc(f)$ is a square in $\Q$.

As the splitting fields considered above are finite extensions of $\Q$ they are number fields. A fundamental object associated with a number field $K$ is its ring of integers,
\[
\mathcal O_K := \{\alpha\in K:\alpha \text{ is a root of a monic polynomial in }\Z[x]\}.
\]
Thus, this ring consists of the elements of $K$ that are integral over $\Z$, and it plays in $K$ a similar role to that of $\Z$ inside $\Q$. In general, $\mathcal O_K$ need not have unique factorization of elements into irreducibles. However, unique factorization is recovered at the level of ideals: for every rational prime $p$, the ideal $p \mathcal O_K$ admits a unique factorization
\[
  p \mathcal O_K = \mathfrak p_1^{e_1}\cdots\mathfrak p_g^{e_g}
\]
into prime ideals of \(\mathcal O_K\). Then, we say that $p$ \emph{ramifies} in $K$ if $e_i>1$ for some $i$. If $[K:\Q]=n$ and $\omega_1,\dots,\omega_n$ is an integral basis of $\mathcal O_K$, we define the discriminant of $K$ as
\[
  \Disc(K):=
  \det\!\left(
    \operatorname{Tr}_{K/\Q}(\omega_i\omega_j)
  \right)_{1\leq i,j\leq n}.
\]
This is an integer independent of the choice of integral basis which detects the ramified rational primes: $p$ ramifies in $K$ if and only if $p \mid \Disc(K)$.
Moreover, Minkowski's discriminant bound implies that $\vert \Disc(K) \vert >1$ whenever $K \neq \Q$. Consequently, every nontrivial number field is ramified at
some finite prime.

The field and polynomial discriminants are related in the following way. If $f(x)\in\Z[x]$ is monic and separable, then every finite prime ramifying in $\Q_{f(x)}$ divides $\Disc(f(x))$. Hence the prime divisors of the polynomial discriminant of $f(x)$ provide a convenient finite set containing all finite primes at which the splitting field of $f(x)$ may ramify.

Let $L_1,\dots,L_s$ be finite extensions of $\Q$. We say that they are jointly linearly disjoint over $\Q$ if
\[
[L_1\cdots L_s:\Q]=\prod_{i=1}^s [L_i:\Q],
\]
where $L_1\cdots L_s$ denotes the compositum of the extensions $L_1,\dots,L_s$. In particular, for two finite Galois extensions $L_1,L_2$ of $\Q$, linear disjointness is equivalent to $L_1\cap L_2 = \Q $. In this case,
\[
  \Gal(L_1L_2/\Q)  \cong \Gal(L_1/\Q)\times\Gal(L_2/\Q).
\]

We will also use the standard fact that a finite prime unramified in each of several finite extensions of $\Q$ remains unramified in their compositum; equivalently, a finite prime ramifying in a compositum must ramify in at least one of the constituent fields. Together with the observations above, this yields the following consequence, which will be used frequently throughout the paper.

\begin{lemma}\label{lem:ramification}
Let $L_1,\dots,L_s$ be finite Galois extensions of $\Q$, and assume that no finite prime ramifies in two distinct fields $L_i$ and $L_j$. Then the fields are jointly linearly disjoint over $\Q$, and
\[
  \Gal(L_1\cdots L_s/\Q)  \cong \prod_{i=1}^s\Gal(L_i/\Q).
\]
\end{lemma}

\begin{proof}
Let $M=L_1\cdots L_{s-1}$ and $E=M\cap L_s$. Since $M/\Q$ and $L_s/\Q$ are Galois, so is $E/\Q$.  As a finite prime can ramify in $M$ only if it ramifies in at least one of $L_1,\dots,L_{s-1}$, by hypothesis, no finite prime ramifies in both $M$ and $L_s$. Since any finite prime ramifying in $E$ would necessarily ramify in both $M$ and $L_s$, the extension $E/\Q$ is unramified at every finite prime. Minkowski's discriminant bound then implies that $E=\Q$. Thus, $M$ and $L_s$ are linearly disjoint over $\Q$. The result now follows by induction.
\end{proof}

\begin{corollary}\label{cor:coprime-disc}
Let $f_1(x),\dots,f_s(x)\in\Z[x]$ be monic separable polynomials. If the polynomial discriminants $\Disc(f_i(x))$ are pairwise coprime, then the splitting fields $\Q_{f_1(x)},\dots, \Q_{f_s(x)}$ are jointly linearly disjoint over $\Q$.
\end{corollary}

\begin{proof}
Every finite prime ramifying in $\Q_{f_i(x)}$ divides $\Disc(f_i(x))$. Since the discriminants $\Disc(f_i(x))$ are pairwise coprime, no finite prime ramifies in two distinct splitting fields. Then, the result follows from Lemma~\ref{lem:ramification}.
\end{proof}

We finish this section with an elementary prime-avoidance result that will allow us to produce splitting fields with disjoint ramification from families of polynomials whose discriminants depend polynomially on an integer parameter.

\begin{lemma}\label{lem:prime-avoidance}
Let $M$ be a nonzero integer and let $f(t)\in\Z[t]$ be a nonconstant polynomial. If for every prime $p$ there exists an integer $a_p$ such that $p\nmid f(a_p)$, then there exist infinitely many positive integers $t_1,t_2,\dots$ such that $\vert f(t_i) \vert >1$, $\gcd(f(t_i),M)=1$, and the integers $f(t_i)$ are pairwise coprime.
\end{lemma}

\begin{proof}
Suppose that $t_1,\dots,t_k$ have been chosen and set
\[
N= \vert M \vert \prod_{i=1}^k \vert f(t_i) \vert.
\]
For each prime $p\mid N$, choose an integer $a_p$ such that $p\nmid f(a_p)$. By the Chinese remainder theorem, there exists an integer $a$ satisfying
$a\equiv a_p\pmod p$ for all $p\mid N$. Then every integer $t$ satisfying 
\begin{equation}\label{progresion}
t\equiv a\pmod{ \prod_{p\mid N}p }
\end{equation} 
also satisfies $t\equiv a_p\pmod p$ for every prime $p\mid N$. Hence $f(t)\equiv f(a_p)\not\equiv 0\pmod p$
for every prime $p\mid N$, and therefore $ \gcd(f(t),N)=1$.

Since $f(x)$ is nonconstant, $\vert f(t) \vert$ is unbounded on the arithmetic progression (\ref{progresion}). Thus, we may choose a positive integer $t$ in this progression such that $\vert f(t) \vert >1$. Setting $t_{k+1}=t$, we obtain $ \gcd(f(t_{k+1}),M)=1$ and $\gcd(f(t_{k+1}),f(t_i))=1$ for all $1\le i\le k$. Finally, iterating this construction we obtain the result.
\end{proof}

\begin{remark}
The same argument can be applied to finitely many families of polynomials simultaneously, under the analogous prime-avoidance hypotheses, while avoiding any prescribed finite set of primes. We will use this observation without further comment.
\end{remark}

Graph theory provides a rich source of polynomials, several of which encode important combinatorial information about the underlying graph and provide
natural objects to which the preceding algebraic tools may be applied. In the next sections, we focus on two important examples: the chromatic polynomial and the domination polynomial.

%%%%%%%%%%%%%%%%%%%%%%%%%%%%%

\section{Galois groups of chromatic polynomials}\label{sec:3}

Let $G=(V,E)$ be a finite simple graph with vertex set $V$ and edge set $E$. A proper $q$-coloring of $G$ is a map $\chi : V \to \{1,\dots,q\}$ such that adjacent vertices receive different colors. In order to count such colorings, Birkhoff \cite{Birkhoff1912} and Whitney \cite{Whitney1932} introduced the chromatic polynomial $P(G,x) \in \Z[x]$ of $G$. It is characterized by the property that, for every positive integer $q$, the value $P(G,q)$ is the number of proper $q$-colorings of $G$.
In particular, $P(G,x)$ is monic of degree $\vert V \vert$.

\subsection{Some chromatic realizations}\label{sec:3.1}

The first systematic study of Galois groups of chromatic polynomials was carried out by Morgan in \cite{Morgan2010, Morgan2012}. Among the usual families of graphs, cycles provide some of the simplest nontrivial examples. Indeed, the chromatic polynomial of the cycle $C_r$ is 
\[
  P(C_r,x)=(x-1)^r+(-1)^r(x-1),
\]
and the non-rational roots generate the cyclotomic field $\Q(\zeta_{r-1})$. Hence
\begin{equation}\label{eq:cycle-galois}
  \Gal(P(C_r,x))\cong(\Z/(r-1)\Z)^\times.
\end{equation}
This provides a convenient source of explicit abelian chromatic Galois groups. For instance, taking $r=p+1$, where $p$ is prime, yields $\Gal(P(C_r,x))\cong (\mathbb Z/p\mathbb Z)^\times\cong \mathbb Z/(p-1)\mathbb Z$, and hence realizes the cyclic group of order $p-1$. However, this construction cannot produce any nontrivial group of odd order, since $\vert (\mathbb Z/n\mathbb Z)^\times \vert =\varphi(n)$ is even for every $n>2$.

The realization of groups of odd order is considerably more subtle. In \cite{Morgan2012}, Morgan introduced a family of graphs $G_{p,q,r}$ with chromatic polynomial
\begin{equation}\label{eq:morgan-general}
 P(G_{p,q,r},x) = P(K_r,x) g_{p,q,r}(x)
 =P(K_r,x) ((x-r)^p+(r-q-1)(x-r+1)^{p-1}).
\end{equation}
As the complete-graph factor has chromatic polynomial $P(K_r,x) = x(x-1) \cdots (x-r+1)$, only the factor $g_{p,q,r}(x)$ contributes to the splitting field, and 
$\Gal(P(G_{p,q,r},x)) = \Gal(g_{p,q,r}(x))$. Then, taking $p=3$ and $r=c^2+c+8+q$, $c \in \N$, the resulting cubic factor $g_{3,q,c^2+c+8+q}(x)$ is irreducible and
\begin{equation}\label{eq:cycle-galois}
  \Gal(P(G_{3,q,c^2+c+8+q},x)) \cong A_3\cong \Z /3\Z.
\end{equation}
Using a related construction, Gilmore found a graph whose chromatic polynomial has Galois group isomorphic to $\mathbb Z/5\mathbb Z$ \cite[Proposition 32]{DelbourgoMorgan2019}. Finding a graph whose chromatic polynomial has Galois group $\mathbb Z/n\mathbb Z$ for odd $n>5$ remains an open problem \cite[Problem 2]{Morgan2012}.

A rather different picture emerges for $\theta$-graphs, whose chromatic polynomials were also studied computationally by Morgan. Recall that the $\theta$-graph $\theta_{a,b,c}$ can be constructed by identifying a path of length $c$ in $C_{a+c}$ with a path of length $c$ in $C_{b+c}$. Among the $\theta$-graphs of order at most $19$, Morgan found that all but two have a chromatic polynomial with an irreducible factor whose Galois group is symmetric \cite[Appendix B]{Morgan2012}.
Beyond these computations, Delbourgo and Morgan proved that the equal-length theta graph $\theta_{a,a,a}$ has Galois group 
\begin{equation}\label{eq:theta-galois}
\Gal(P(\theta_{a,a,a},x)) \cong S_{3(a-1)}
\end{equation}
for every $a\geq 2$ \cite[Theorem 1.3]{DelbourgoMorgan2014}.

\subsection{Products via $K_1$-gluings}\label{sec:3.2}

All the examples considered above are connected graphs. It is therefore natural to ask how the Galois group behaves when disconnected graphs are allowed. Let
$G=G_1\sqcup\cdots\sqcup G_r$ be the decomposition of $G$ into connected components. Since the components may be colored independently, we have 
\begin{equation}\label{eq:chromatic-components}
 P(G,x)=\prod_{i=1}^rP(G_i,x).
\end{equation}
Consequently, the splitting field of $P(G,x)$ is the compositum of the splitting fields of the polynomials $P(G_i,x)$ attached to the connected components. Thus, disconnected graphs provide a natural way of combining chromatic Galois groups. However, in the chromatic setting, the same realizations can be obtained by connected graphs using $K_1$-gluings. Recall that if two graphs $G_1$ and $G_2$ each contain a copy of the complete graph $K_r$, their $K_r$-gluing is obtained by identifying these two copies of $K_r$. Its chromatic polynomial satisfies
\begin{equation}\label{eq:chromatic-clique-sum}
P(G_1\cup_{K_r}G_2,x) = \frac{P(G_1,x)P(G_2,x)}{P(K_r,x)},
\end{equation}
and as $P(K_r,x)$ splits completely over $\mathbb Q$, this operation does not change the compositum of the splitting fields contributed by $G_1$ and $G_2$.
 
\begin{proposition}\label{prop:chromatic-connected}
Every chromatic Galois realization admits a connected realization. 
\end{proposition}

\begin{proof}
Let $G=G_1\sqcup\cdots\sqcup G_s$ be the decomposition of a realizing graph into connected components. Choose one vertex in each component and form a connected graph $G'$ by taking successive $K_1$-gluings. Repeated use of \eqref{eq:chromatic-clique-sum} gives
\[
 P(G',x)=\frac{\prod_{i=1}^sP(G_i,x)}{x^{s-1}}
       =\frac{P(G,x)}{x^{s-1}}.
\]
Then, as the removed factor splits over $\Q$, $P(G,x)$ and $P(G',x)$ have the same splitting field and consequently the same Galois group.
\end{proof}

Now, we return to Morgan's cubic family introduced above. Combined with the preceding result and a suitable choice of ramification, we can prove that this family yields connected realizations of elementary abelian $3$-groups of arbitrary rank.

\begin{proposition}\label{pro:realizationC3}
For every $m\geq 0$, there exists a connected graph $G$ such that
\[
\Gal(P(G,x))\cong (\Z/3\Z)^m.
\]
\end{proposition}

\begin{proof}
The case $m=0$ is trivial: as $P(K_r,x)$ splits completely over $\Q$, this realizes the trivial group as the Galois group of a chromatic polynomial. Then, we may assume that $m>0$. Let
\[
g_{p,q,r}(x) =(x-r)^p+(r-q-1)(x-r+1)^{p-1}
\]
the nontrivial factor in \eqref{eq:morgan-general} and recall that $\Gal(P(G_{p,q,r},x)) = \Gal(g_{p,q,r}(x))$. Taking $p=3$, $q=0$, and $r=c^2+c+8$, the discriminant of $g_{3,0,c^2+c+8}(x)$ is
\[
f_3(c)^2 = \Disc(g_{3,0,c^2+c+8}(x))
=
\bigl((c^2+c+7)(2c+1)\bigr)^2,
\]
and the polynomial $f_3(c)$ has no fixed prime divisor; indeed, for $p \neq 3$ we consider $f_3(1)=27$ and for $p=3$ we consider $f_3(2)=65$.
Then, applying Lemma \ref{lem:prime-avoidance} with $M = 1$, we may choose positive integers $c_1,\ldots,c_m$ such that the values $f_3(c_i)$ are pairwise coprime. Consequently, the polynomial discriminants of the cubic irreducible polynomials $g_{3,0,c_i^2+c_i+8}(x)$ are pairwise coprime and by Corollary \ref{cor:coprime-disc}, their splitting fields are jointly linearly disjoint over $\Q$. Thus, as $\Gal(g_{3,0,c_i^2+c_i+8}(x)) \cong \Z/3\Z$, the compositum $L$ of all these splitting fields satisfies $\Gal(L/\Q) \cong (\Z/3\Z)^m$.

Now, consider the graph
\[
X=G_{3,0,c_1^2+c_1+8} \sqcup\cdots\sqcup G_{3,0,c_m^2+c_m+8}.
\]
Since the chromatic polynomial of a disjoint union is the product of the chromatic polynomials of its components \eqref{eq:chromatic-components}, the splitting field of 
\[
P(X,x) = \prod_{i=1}^m P(G_{3,0,c_i^2+c_i+8},x)
\] 
is precisely $L$. Then, $\Gal(P(X,x)) \cong (\Z/3\Z)^m$ and by Proposition \ref{prop:chromatic-connected} there exists a connected graph $G$ with $\Gal(P(G,x)) \cong (\Z/3\Z)^m$, as required.
\end{proof}

More generally, we have the following.

\begin{theorem}\label{thm:chromatic-stability}
Let $\Gamma$ be the Galois group of a chromatic polynomial. Then, for every $n,m\geq0$, there exists a connected graph $G$ such that
\[
\Gal(P(G,x)) \cong \Gamma\times(\Z/2\Z)^n\times(\Z/3\Z)^m.
\]
\end{theorem}

\begin{proof}
Let $H_1$ be a connected graph such that $\Gamma \cong \Gal(P(H_1,x))$ and let $M$ be the product of the primes dividing the field discriminant $\Disc(\Q_{P(H_1,x)})$, with $M=1$ if there are none. If $n=m=0$ there is nothing to prove, so we may assume that $n+m>0$. 

Note that the family of graphs $G_{p,q,r}$ constructed by Morgan also yields infinitely many realizations of $\mathbb Z/2\mathbb Z$. Indeed, taking $p=2$, $q=0$ and $r \neq 1,5$, the nontrivial factor $g_{2,0,r}(x)$ in \eqref{eq:morgan-general} is an irreducible quadratic polynomial with discriminant
\[
f_2(r) := \Disc(g_{2,0,r}(x))=(r-1)(r-5).
\]
As in the cubic case, the polynomial $f_2(r)$ has no fixed prime divisor: for $p \neq 3$ we use that $f_2(2)=-3$ and for $p=3$ we use that $f_2(3)=-4$. Hence, by Lemma \ref{lem:prime-avoidance}, we may choose positive integers $ r_1,\ldots,r_n$, such that the integers $f_2(r_i)$ are pairwise coprime and prime to $M$. Then, by Corollary \ref{cor:coprime-disc} the splitting fields $\Q_{g_{2,0,r_1}(x)}, \ldots, \Q_{g_{2,0,r_n}(x)}$ are jointly linearly disjoint over $\Q$.  Therefore, since  $\Gal(g_{2,0,r_i}(x))= \Gal(P(G_{2,0,r_i},x))\cong \Z/2\Z$, considering the graph
\[
Y=G_{2,0,r_1} \sqcup\cdots\sqcup G_{2,0,r_n},
\]
we can construct a connected graph $H_2$ such that $\Gal(P(H_2,x)) \cong (\Z/2\Z)^n$ and no prime dividing $M$ ramifies in the splitting fields $\Q_{P(H_2,x)}$. 

On the other hand, the construction in Proposition \ref{pro:realizationC3}, carried out there with $M=1$, works equally well for an arbitrary positive integer $M$. Hence there exists a connected graph $H_3$ with $\Gal(P(H_3,x))\cong (\Z/3\Z)^m$ and $\Q_{P(H_3,x)}$ unramified at every prime dividing $M$.

Finally, set $L_1=\Q_{P(H_1,x)}$, $L_2=\Q_{P(H_2,x)}$ and $L_3=\Q_{P(H_3,x)}$. Since $[L_2:\Q]=2^n$ and $[L_3:\Q]=3^m$, the degree $ [L_2\cap L_3:\Q]$ divides both $2^n$ and $3^m$. Hence $L_2\cap L_3=\Q$, then $L_2$ and $L_3$ are linearly disjoint over $\Q$. Moreover, both $L_2$ and $L_3$ are unramified at every prime dividing $M$, and hence so is their compositum $L_2L_3$. As every finite prime ramifying in $L_1$ divides $M$, then $L_1$ and $L_2L_3$ have disjoint finite ramification, and therefore they are linearly disjoint over $\Q$. Consequently, $\Gal(L_1L_2L_3/\Q) \cong \Gamma\times(\Z/2\Z)^n\times(\Z/3\Z)^m
$.

Now consider the disjoint union 
\[
Z=H_1\sqcup H_2\sqcup H_3.
\] 
Its splitting field is $L_1L_2L_3$, then $\Gal(P(Z,x)) \cong \Gamma\times(\Z/2\Z)^n\times(\Z/3\Z)^m$. Thus, by Theorem \ref{prop:chromatic-connected}, this realization admits a connected realization. Hence there exists a connected graph $G$ such that
\[
\Gal(P(G,x)) \cong \Gamma\times(\Z/2\Z)^n\times(\Z/3\Z)^m,
\]
as required.
\end{proof}

This theorem may be applied to any known chromatic Galois realization. In particular, using the the Galois groups of equal-length theta graphs obtained by Delbourgo--Morgan, we have the following result.

\begin{corollary}\label{cor:chromatic-S3k}
For every $n,m\geq0$ and $k\geq1$, there exists a connected graph $G$ such that 
\[
\Gal(P(G,x)) \cong (\Z/2\Z)^n\times(\Z/3\Z)^m \times S_{3k}.
\]
\end{corollary}

%%%%%%%%%%%%%%%%%%%%%%%%%%%%%%%%%%%%%%%%%%%%%%%%%%%%%%%%%%%%

\section{Galois groups of domination polynomials}\label{sec:4}

Let $G=(V,E)$ be a finite simple graph. A subset $S\subseteq V$ is dominating if every vertex in $V\setminus S$ has a neighbor in $S$. Let $d_k(G)$ denote the number of dominating sets of cardinality $k$, and let $\gamma(G)$ be the domination number of $G$, i.e., the minimum cardinality of a dominating set. The domination polynomial of $G$ is defined by
\[
 D(G,x)=\sum_{k=\gamma(G)}^{|V|}d_k(G)x^k.
\]
As usual, rational linear factors may be omitted when considering splitting fields or Galois groups. In particular, we will omit the factor $x^{\gamma(G)}$.

\subsection{First domination realizations}\label{sec:4:1}

We begin with two basic families for which the domination polynomial admits an explicit description: complete graphs and stars.
For a complete graph $K_n$, every nonempty vertex subset is dominating; hence
\[
D(K_n,x)=(1+x)^n-1.
\]
Its roots are $\zeta_n^j-1$, $0\leq j\leq n-1$, where $\zeta_n$ is a primitive $n$th root of unity. Therefore, the splitting field is $\Q(\zeta_n)$, and we obtain the analogue of the corresponding result for cycles in the chromatic case.

\begin{proposition}\label{prop:complete}
For every $n\ge2$,
\[
  \Gal(D(K_n,x))\cong(\Z/n\Z)^\times.
\]
In particular, if $p$ is prime, then $K_p$ realizes $C_{p-1}$.
\end{proposition}

Now, let $K_{1,n}$ be the star with $n$ leaves. A dominating set either contains the center, in which case the leaves may be chosen freely, or omits the center, in which case every leaf must be chosen. Hence
\[
 D(K_{1,n},x)=x(1+x)^n+x^n
 =x\bigl((1+x)^n+x^{n-1}\bigr).
\]
The Galois groups of star graphs already exhibit a markedly different behavior.

\begin{theorem}\label{thm:stars}
For every $n\ge2$,
\[
 \Gal(D(K_{1,n},x))\cong
 \begin{cases}
  S_n,& n\not\equiv5\pmod6,\\[2mm]
  \Z / 2 \Z \times S_{n-2},& n\equiv5\pmod6.
 \end{cases}
\]
\end{theorem}

\begin{proof}
For a nonzero root, consider the change of variables $z=\frac{1+x}{x}$. Then $(1+x)^n+x^{n-1}=0$ if and only if
\[
g_n(z):=z^n+z-1=0.
\]
Since the transformation is rational and invertible over $\Q$, with inverse $x=\frac{1}{z-1}$, the splitting field of $D(K_{1,n},x)$ coincides with the splitting field of $g_n(z)$. It therefore suffices to compute the Galois group of $g_n(z)$.

If $n$ is even, the substitution $z=-y$ transforms $g_n(z)$ into $y^n-y-1$. By a result of Osada~\cite[Corollary 3]{Osada1987} this polynomial has Galois group $S_n$.

Assume now that $n$ is odd. The same substitution transforms $g_n(z)$, up to multiplication by $-1$, into $y^n+y+1$. In ~\cite[Proposition~2.3]{DelbourgoMorgan2014}, Delbourgo--Morgan compute the Galois group of $f_a(y)=y^{a+1}+y+1=y^n+y+1$ 
\[
 \Gal(f_a(y))\cong
 \begin{cases}
  S_{a+1}=S_n,&a\not\equiv1\pmod3,\\
  \Z/2 \Z \times S_{a-1}= \Z/2 \Z \times S_{n-2},&a\equiv1\pmod3.
 \end{cases}
\]
Since $n$ is odd, the congruence $a\equiv1\pmod3$ is equivalent to $n\equiv5\pmod6$. 
\end{proof}

The preceding proof highlights an interesting link with the chromatic case. The polynomial $y^{a+1}+y+1$, used by Delbourgo--Morgan in obtaining chromatic realizations of $S_{3k}$ from equal-length theta graphs, also appears here naturally from the domination polynomials of stars, yielding realizations of $S_n$ for many values of $n$.

The star family is not isolated. Let $G$ be a graph with vertex set $V=\{v_1,\ldots,v_r\}$, and let $m_1,\ldots,m_r$ be positive integers. Following the notation of Oboudi~\cite{Oboudi2016}, write
\[
 G(v_1^{m_1},\ldots,v_r^{m_r})
\]
for the sunlike graph obtained by attaching $m_i$ new leaves to $v_i$ for each $i$. Oboudi's factorization theorem~\cite[Theorem~3]{Oboudi2016} states that
\begin{equation}\label{eq:sunlike}
 D(G(v_1^{m_1},\ldots,v_r^{m_r}),x)
 =\prod_{i=1}^rD(K_{1,m_i},x).
\end{equation}
A notable feature of this identity is that the edges of the base graph $G$ do not appear on the right-hand side.

\begin{corollary}\label{cor:sunlike-composita}
Let $H=G(v_1^{m_1},\ldots,v_r^{m_r})$, and $L_i$ denote the splitting field of $D(K_{1,m_i},x)$. Then $\Q_{D(H,x)}=L_1\cdots L_r$. Moreover, if the fields $L_1,\ldots,L_r$ are jointly linearly disjoint, then
\[
 \Gal(D(H,x))\cong
 \prod_{i=1}^r\Gal(D(K_{1,m_i},x)).
\]
Moreover, if $G$ is connected, then $H$ is connected, and if $G$ is a tree, then $H$ is a tree.
\end{corollary}

\begin{proof}
The splitting-field statement follows immediately from \eqref{eq:sunlike}, while the direct-product assertion follows from the joint linear disjointness of the fields $L_1,\ldots,L_r$. Attaching pendant vertices preserves connectedness and, when the base graph is a tree, preserves acyclicity.
\end{proof}

The uniform case gives an especially clean propagation of the star theorem.

\begin{corollary}\label{cor:star-groups-trees}
Let $m\ge2$ and let $T$ be any tree with vertex set $V=\{v_1,\ldots,v_r\}$. If $T^{(m)}=T(v_1^m,\ldots,v_r^m)$,
then
\[
 \Gal(D(T^{(m)},x))\cong
 \begin{cases}
 S_m,&m\not\equiv5\pmod6,\\[2mm]
 \Z/2\Z \times S_{m-2},&m\equiv5\pmod6.
 \end{cases}
\]
Consequently, every Galois group realized by a star in Theorem~\ref{thm:stars} is realized by infinitely many pairwise nonisomorphic trees. In fact, these realizations may be chosen to be caterpillars.
\end{corollary}

\begin{proof}
By \eqref{eq:sunlike},
\[
 D(T^{(m)},x)=D(K_{1,m},x)^r,
\]
Then, $D(T^{(m)},x)$ and $D(K_{1,m},x)$ have the same splitting field, and Theorem~\ref{thm:stars} applies. Taking $T=P_r$, the path graph on $r$ vertices, for $r=1,2,\ldots$, gives caterpillars of strictly increasing order and hence infinitely many pairwise nonisomorphic tree realizations.
\end{proof}

\subsection{Connected products for domination polynomials}\label{sec:4.2}

Let $G=G_1\sqcup\cdots\sqcup G_r$ be the decomposition of $G$ into connected components. Since a set dominates $G$ exactly when its intersection with each connected component dominates that component, the domination polynomial is multiplicative over connected components:
\begin{equation}\label{eq:domination-components}
D(G,x)=\prod_{i=1}^r D(G_i,x).
\end{equation}

Thus, as in the chromatic case, the splitting field of $D(G,x)$ is the compositum of the splitting fields of the polynomials $D(G_i,x)$. However, unlike the chromatic case, there is no known general procedure for replacing an arbitrary disconnected domination realization by a connected one while preserving its Galois group. In particular, we do not have an analogue of the $K_1$-gluing argument of Proposition~\ref{prop:chromatic-connected}. Nevertheless, under additional hypotheses on the components, there are constructions that preserve the product while producing a connected graph.

For example, Corollary~\ref{cor:sunlike-composita} provides a mechanism for producing connected realizations of products of star Galois groups. More generally, the factorization \eqref{eq:sunlike} underlying this corollary is a special case of a construction of Akbari, Alikhani, and Peng~\cite[Theorem~7]{AkbariAlikhaniPeng2010}.

Let $H_1,\ldots,H_r$ be graphs with chosen universal vertices $u_1,\ldots,u_r$, where a vertex $u$ of a graph $H$ is universal if it is adjacent to every other vertex of $H$. Identifying each $u_i$ with the corresponding vertex $v_i$ of an arbitrary base graph $G$ on $r$ vertices, and denoting the resulting graph by $H$, we obtain
\begin{equation}\label{eq:productAAP}
D(H,x)=\prod_{i=1}^r D(H_i,x).
\end{equation}
Thus, under the universal-vertex hypothesis, the components can be connected through an arbitrary base graph while preserving the product of their domination polynomials. The sunlike construction is recovered by taking $H_i=K_{1,m_i}$ and $u_i$ to be the center of the star.

The construction above already yields connected Galois realizations that cannot be obtained from complete graphs. To illustrate this, first note that complete graphs cannot realize elementary abelian $2$-groups of arbitrary rank.
Let $r \geq 1$ and suppose that $\Gal(D(K_n,x)) \cong (\Z/n\Z)^\times\cong(\Z/2\Z)^r$ for some complete graph $K_n$. Write 
\[ 
n=2^a\prod_{i=1}^s p_i^{e_i},
\]
where the $p_i$ are distinct odd primes. For an odd prime $p$, the group $(\Z/p^e\Z)^\times$ is cyclic of order $p^{e-1}(p-1)$. If it is to be an elementary abelian $2$-group, we must have $e=1$ and $p-1\leq2$, then $p=3$. Hence the odd part of $N$ is either $1$ or $3$. On the other hand, $(\Z/2^a\Z)^\times$ is elementary abelian only for $a\leq3$, since for $a\geq4$ it contains an element of order $2^{a-2}$. Therefore $n \mid 24$, and consequently $r\leq3$. Conversely,
\[
\Gal(D(K_{24},x)) \cong (\Z/24\Z)^\times\cong(\Z/2\Z)^3,
\]
then $K_{24}$ realizes the rank three. Rank one and two are realized by $K_3$ and $K_8$ respectively.

\begin{proposition}\label{prop:C2four-dom}
There exists a connected graph whose domination polynomial has Galois group $(\Z/2\Z)^4$, and no complete graph has this Galois group.
\end{proposition}

\begin{proof}
Take $H_1=K_{24}$, $H_2=K_{1,2}$, and let the base graph be $K_2$. Choose any vertex of $K_{24}$ and the center of $K_{1,2}$ as the corresponding universal vertices. By \eqref{eq:productAAP}, the resulting graph $H$ is connected and
\[
D(H,x)=D(K_{24},x)D(K_{1,2},x).
\]
Now, as $D(K_{1,2},x)=x(x^2+3x+1)$, its splitting field is $\Q(\sqrt5)$ and $\Gal(D(K_{1,2},x)) \cong \Z/2\Z$.  On the other hand, from \S\ref{sec:4:1} and the previous discussion, the splitting field of $D(K_{24},x)$ is the cyclotomic field $\Q(\zeta_{24})$ and $\Gal(D(K_{24},x))\cong(\Z/2\Z)^3$. Since $\Q(\sqrt5)$ is ramified only at $5$ whereas $\Q(\zeta_{24})$ is ramified only at $2$ and $3$, these fields are linearly disjoint. Hence 
\[
\Gal(D(H,x)) \cong \Gal(D(K_{24},x)) \times \Gal(D(P_3,x)) \cong (\Z/2\Z)^4.
\]
The last assertion follows from the preceding observation.
\end{proof}

This raises a natural problem.

\begin{question}\label{q:C2powers-dom}
Is $(\Z/2\Z)^r$ realizable as the Galois group of a domination polynomial for every $r\geq1$? Can it always be realized by a connected graph?
\end{question}
The product construction reduces this question to finding sufficiently many connected $(\Z/2\Z)$-realizations with universal vertices whose quadratic splitting fields are linearly disjoint. This contrasts with the chromatic setting, where $(\Z/2\Z)^r$ admits connected realizations for every $r\geq0$.

This method applies to a considerably larger class of factors, including complete graphs (then, we can also prove that every Galois group realized by a complete graph in Proposition~\ref{prop:complete} is realized by infinitely many pairwise nonisomorphic graphs). However, it does not provide a way of combining arbitrary domination polynomials, since the existence of a universal vertex is essential. A more flexible connected-product construction is obtained from the theory of irrelevant edges and domination-covered vertices. 

A vertex $v$ of $G$ is called domination-covered if every dominating set of $G-v$ contains a neighbor of $v$. Equivalently, there exists $u \in N_G(v)$ such that $N_G[u] \subseteq N_G[v]$. An edge $e$ is irrelevant when $D(G,x)=D(G-e,x)$. Kotek--Preen--Simon--Tittmann--Trinks prove that an edge $uv$ is irrelevant precisely when both endpoints are domination-covered in $G-uv$~\cite[Theorem 3.11]{KotekEtAl2012}. Then, for graphs possessing the suitable vertices, this result gives the following connected product principle.

\begin{proposition}\label{prop:connected-product}
Let $G_1$ and $G_2$ be connected graphs with domination-covered vertices $u_1$ and $u_2$ respectively. Let $H$ be obtained from $G_1\sqcup G_2$ by adding the edge $u_1u_2$. Then $H$ is connected and
\[
  D(H,x)=D(G_1,x)D(G_2,x).
\]
Moreover, if the splitting fields $L_i=\Q_{D(G_i,x)}$ satisfy $L_1\cap L_2=\Q$, then
\[
  \Gal(D(H,x))\cong
  \Gal(D(G_1,x))\times\Gal(D(G_2,x)).
\]
\end{proposition}
\begin{proof}
Deleting the new edge from $H$ we obtain $G_1\sqcup G_2$, in which both $u_1$ and $u_2$ are domination-covered. Hence the new edge $u_1u_2$ is irrelevant and he polynomial identity follows from \eqref{eq:domination-components}. The splitting field of $D(H,x)$ is $L_1L_2$, and the Galois-group statement follows from linear disjointness of $L_1$ and $L_2$.
\end{proof}

The construction can be iterated. Indeed, chosen domination-covered vertices of finitely many graphs may be connected according to the edges of a tree. The neighborhood inclusions witnessing domination-coveredness remain valid as new incident edges are added, since the closed neighborhood of each chosen vertex can only increase. Thus every added edge is irrelevant at the moment it is introduced. Consequently, if the corresponding splitting fields are jointly linearly disjoint, the resulting connected graph realizes the direct product of their Galois groups.

\begin{remark}\label{rmk:D5}
This construction is genuinely more flexible than the universal-vertex factorization at the level of graph realizations.
Let $F$ be the graph obtained from $K_4$ by choosing a vertex $v$ and adjoining three new vertices $w_1,w_2,w_3$ together with the edges $vw_1, w_1w_2, w_2w_3$. A direct computation gives
\begin{equation}\label{eq:dihedral5}
D(F,x)=x^2(x^5+7x^4+20x^3+29x^2+21x+5).
\end{equation}
Using MAGMA, we can verify that the quintic factor is irreducible with discriminant $47^2$ and Galois group the dihedral group $D_5$ of order 10. The graph $F$ has no universal vertex, but $w_2$ is domination-covered, since $N_F[w_3]\subseteq N_F[w_2]$. Hence $F$ cannot be used as a factor in the universal-vertex construction, whereas it can be used in Proposition \ref{prop:connected-product}. For example, every vertex of $K_3$ is domination-covered, and $\Q_{(D(K_3,x))}=\Q(\zeta_3)$. The splitting field of $D(F,x)$ is the splitting field of the quintic factor above and can ramify only at $47$, whereas $\Q(\zeta_3)$ ramifies only at $3$. Hence the two fields are linearly disjoint. Joining $w_2$ to a vertex of $K_3$ by an irrelevant edge therefore gives a connected graph $G$ satisfying
$\Gal(D(G,x))\cong D_5\times \Z/2\Z$ (Figure \ref{fig:connected-D5}).
%This construction is illustrated in Figure \ref{fig:connected-D5}.

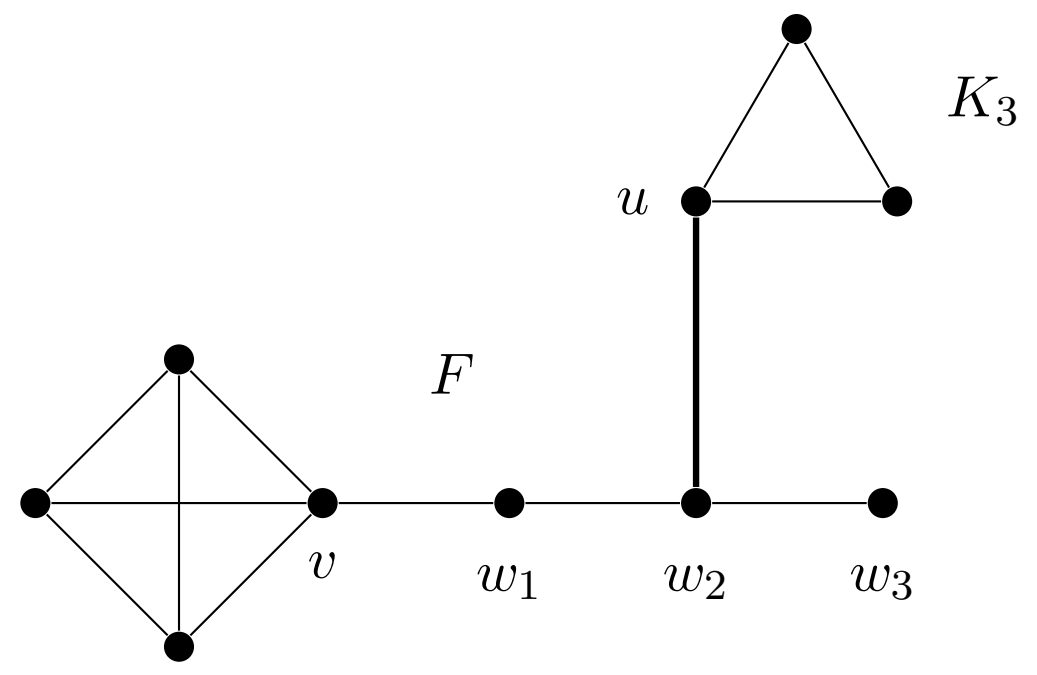
\begin{figure}[ht]
\centering
\begin{tikzpicture}[
  scale=0.95,
  every node/.style={circle,fill=black,inner sep=2pt}
]
% The graph F: a K_4 with an attached path
\node (a) at (-3,0) {};
\node (b) at (-2,1) {};
\node (c) at (-2,-1) {};
\node (v) at (-1,0) {};

\node (w1) at (0.3,0) {};
\node (w2) at (1.6,0) {};
\node (w3) at (2.9,0) {};

\draw (a)--(b)--(v)--(c)--(a);
\draw (a)--(v);
\draw (b)--(c);
\draw (v)--(w1)--(w2)--(w3);

% The triangle K_3
\node (u) at (1.6,2.1) {};
\node (p) at (3,2.1) {};
\node (q) at (2.3,3.3) {};

\draw (u)--(p)--(q)--(u)--cycle;

% The irrelevant edge
\draw[very thick] (w2)--(u);

% Labels
\node[fill=none,below=5pt] at (v) {$v$};
\node[fill=none,below=5pt] at (w1) {$w_1$};
\node[fill=none,below=5pt] at (w2) {$w_2$};
\node[fill=none,below=5pt] at (w3) {$w_3$};

\node[fill=none,left=5pt] at (u) {$u$};
\node[fill=none] at (-0.1,0.9) {$F$};
\node[fill=none] at (3.6,2.8) {$K_3$};

\end{tikzpicture}
\caption{The connected graph $G$ with irrelevant edge $uw_2$.}
\label{fig:connected-D5}
\end{figure}

\end{remark}

\begin{remark}
Although the domination-covered condition is substantially weaker than the existence of a universal vertex, it is still far from ubiquitous. 
No cycle $C_s$ with $s\ge4$ has a domination-covered vertex: after deleting a vertex $v$, the internal vertices of the resulting path form a dominating set containing neither neighbor of $v$. Similarly, no complete bipartite graph $K_{a,b}$ with $a,b\ge2$ has a domination-covered vertex; after deleting a vertex in one part, the remaining vertices of that part form a dominating set containing no neighbor of the deleted vertex.
Thus, although the irrelevant-edge construction is more flexible than the universal-vertex construction, it still does not provide a general method for connecting arbitrary components while preserving the product of their domination polynomials. 
\end{remark}

Consequently, connected realizability remains a genuine structural issue rather than a formal consequence of multiplicativity.

\begin{question}\label{q:connected-classes}
Does every group realizable by a domination polynomial admit a realization by a connected graph?
\end{question}

%%%%%%%%%%%%%%%%%%%%%%%%%%%%%%%%%%%%%%%%%%%

\section{Odd-order groups and Oboudi's conjecture}\label{sec:5}

We now turn from realization results to a possible obstruction. Interestingly, the sunlike graphs considered above arise again in this setting, through Oboudi's conjectural characterization of graphs whose domination polynomial has only real roots.

First, let us introduce a useful general observation that is not specific to graph polynomials.

\begin{lemma}\label{lem:odd-real}
Let $f(x)\in\Q[x]$. If $\Gal(f(x))$ has odd order, then all roots of $f(x)$ are real.
\end{lemma}

\begin{proof}
Since $f(x)$ has rational coefficients, complex conjugation permutes its roots and therefore preserves its splitting field $\Q_{f(x)}$. Its restriction to $\Q_{f(x)}$ thus defines an element $c\in\Gal(f(x))$ whose order is at most two. If $\Gal(f(x))$ has odd order, then $c$ must be the identity. Hence every element of $\Q_{f(x)}$ is fixed by complex conjugation, so $\Q_{f(x)}\subseteq \R$. In particular, every root of $f(x)$ is real.
\end{proof}

Consequently, if a nontrivial group of odd order is realizable as the Galois group of a domination polynomial, then that polynomial must have all its roots real and and nontrivial splitting field. Oboudi's conjecture predicts that this never happens.

The factorization \eqref{eq:sunlike} shows that sunlike graphs provide a natural source of domination polynomials with only real roots. Indeed, when every $m_i\in{1,2}$, all roots are real, since
\[
 D(K_{1,1},x)=x(x+2),\qquad D(K_{1,2},x)=x(x^2+3x+1).
\]
Oboudi conjectured \cite[Conjecture 1]{Oboudi2016} that, apart from isolated vertices, these are the only graphs whose domination polynomial has only real roots.

\begin{conjecture}[Oboudi]\label{conj:oboudi}
If all roots of $D(H,x)$ are real, then
\[
 H=G(v_1^{m_1},\ldots,v_r^{m_r})\sqcup tK_1
\]
for some finite simple graph $G$, some $t\ge0$, and $m_i\in{1,2}$ for every $i$.
\end{conjecture}

Although the conjecture is stated for arbitrary finite simple graphs, its specialization to connected graphs is immediate. If $H$ is connected and has at least two vertices, then necessarily $t=0$ and the base graph $G$ is connected. Thus, in the connected setting, the conjectural conclusion becomes
\[
 H=G(v_1^{m_1},\ldots,v_r^{m_r}),
 \qquad G\text{ connected},\quad m_i\in\{1,2\}.
\]

From the Galois-theoretic point of view, these conjecturally exceptional graphs are especially simple: their domination polynomials have at most a quadratic splitting field.

\begin{proposition}\label{prop:sunlike}
Let
\[
 H=G(v_1^{m_1},\ldots,v_r^{m_r})\sqcup tK_1,
 \qquad m_i\in\{1,2\}.
\]
Then
\[
 \Gal(D(H,x))\cong
 \begin{cases}
  1,&m_i=1\text{ for every }i,\\
  \Z / 2\Z,&m_i=2\text{ for at least one }i.
 \end{cases}
\]
\end{proposition}

\begin{proof}
The isolated vertices contribute only powers of $x$. By \eqref{eq:sunlike}, the remaining factors are $D(K_{1,1},x)=x(x+2)$ and $D(K_{1,2},x)=x(x^2+3x+1)$.
Thus the only possible non-rational roots come from $x^2+3x+1$, whose splitting field is $\Q(\sqrt5)$. If $m_i=1$ for every $i$, all roots are rational and the splitting field is $\Q$. If $m_i=2$ for at least one $i$, then $x^2+3x+1$ occurs as a factor, so the splitting field is $\Q(\sqrt5)$. 
\end{proof}

Combining Lemma \ref{lem:odd-real}, Conjecture \ref{conj:oboudi}, and Proposition \ref{prop:sunlike}, we obtain the following conditional obstruction.

\begin{corollary}\label{cor:odd-oboudi}
Assuming Oboudi's conjecture, the trivial group is the only finite group of odd order realizable as the Galois group of a domination polynomial.
\end{corollary}

Thus, if Oboudi's conjecture is true, the inverse Galois problem for domination polynomials has a negative answer. Conversely, if $\Z/3\Z$, or any nontrivial finite group of odd order, is realized as the Galois group of a domination polynomial, then Lemma \ref{lem:odd-real} shows that the corresponding domination polynomial has only real roots. By Proposition \ref{prop:sunlike}, such a graph cannot belong to the family described in Conjecture \ref{conj:oboudi}, and hence provides a counterexample to Oboudi's conjecture.

The preceding discussion naturally leads to two complementary questions, separating the even- and odd-order cases.

\begin{question}\label{q:even}
Is every finite group of even order realizable as the Galois group of a domination polynomial? Is every such group realizable by a connected graph?
\end{question}

\begin{question}\label{q:odd}
Does any nontrivial finite group of odd order occur as the Galois group of a domination polynomial? In particular, can $\Z/3\Z$ be realized?
\end{question}

\begin{remark}\label{rem:small-computation}
Alikhani--Peng tabulated the domination polynomials of all $143$ connected nonisomorphic graphs of order at most six~\cite{AlikhaniPeng2014}. A computation of the Galois groups of these polynomials using MAGMA produced no example with a nontrivial Galois group of odd order.
\end{remark}

%%%%%%%%%%%%%%%%%%%%%%%%%%%%%%%%%%%%%%%%%%%%%%%%%%%%%%

\section{Inverse Galois problems for small groups}\label{sec:6}

The preceding constructions suggest two complementary ways to organize the inverse Galois problems. We may proceed by the order of the abstract Galois groups, or by the degree of an irreducible factor and hence by the transitive subgroups of a symmetric group.

Let us first organize the inverse problem by the order of the abstract Galois group. For chromatic polynomials, the realizations discussed previously give a complete picture through order six.

\begin{corollary}\label{cor:small-chromatic}
Every finite group of order at most six is the Galois group of the chromatic polynomial of a connected graph.
\end{corollary}

\begin{proof}
As seen in Section \ref{sec:3}, cycles provide realizations of $\Z/2\Z$, $\Z/4\Z$, $\Z/6\Z$, and $(\Z/2\Z)^2$. Morgan and Delbourgo provided connected chromatic realizations of $\Z/3\Z$ and $S_3$~\cite{Morgan2012, DelbourgoMorgan2014}, while Gilmore gave a chromatic realization of $\Z/5\Z$~\cite{DelbourgoMorgan2019, Gilmore2015}. Together with the trivial group, these examples exhaust all isomorphism types of finite groups of order at most six.
\end{proof}

For domination polynomials, the constructions of Section \ref{sec:4} show that complete graphs realize $\Z/2\Z$, $\Z/4\Z$, $(\Z/2\Z)^2$, and $\Z/6\Z$ through $K_3$, $K_5$, $K_8$, and $K_7$, respectively, while $K_{1,3}$ realizes $S_3$.

\begin{corollary}\label{cor:small-domination-even}
Every finite group of even order at most six is the Galois group of the domination polynomial of a connected graph.
\end{corollary}

The only remaining groups of order at most six are $\Z/3\Z$ and $\Z/5\Z$. No realization of either group as the Galois group of a domination polynomial is currently known; assuming Oboudi's conjecture, neither can occur, by Corollary \ref{cor:odd-oboudi}.

We now turn to the second point of view, organizing the problem by the degree of an irreducible factor and hence by the transitive subgroups of the corresponding symmetric group. This is closer to the computational viewpoint of Cameron--Morgan~\cite{CameronMorgan2017}.
If an irreducible factor has degree $d$, its Galois group is a transitive subgroup of $S_d$. Cameron--Morgan found that their computations with chromatic polynomials of rings of cliques realize every transitive permutation group of degree at most four; in degree five only $\Z/5\Z$ was missing, while in degree six they found five of the sixteen transitive groups~\cite[Section~6]{CameronMorgan2017}. Gilmore subsequently produced a biclique whose chromatic polynomial has a quintic factor with Galois group $\Z/5\Z$ and, in fact, examples realizing every transitive subgroup of $S_5$~\cite[Propositions~7.6--7.7]{Gilmore2015}; see also \cite[Proposition~32]{DelbourgoMorgan2019}.

For domination polynomials, we obtain the following partial analogue.

\begin{proposition}\label{prop:small-transitive-dom}
For $d\leq4$, every transitive subgroup of $S_d$ except $\Z/3\Z=A_3$ and $A_4$ occurs as the Galois group of an irreducible factor of the domination polynomial of a connected graph.
\end{proposition}

\begin{proof}
For degree two, $S_2$ occurs for $K_3$. In degree three, $S_3$ occurs for $K_{1,3}$, while $\Z/3\Z$ remains open. The transitive subgroups of $S_4$ are  $\Z/4\Z,V_4 = (\Z/2\Z)^2,D_4,A_4$, and $S_4$. The groups $\Z/4\Z$ and $V_4$ occur as Galois groups of the irreducible quartic factors of $D(K_5,x)$ and $D(K_8,x)$, respectively, while $S_4$ occurs for $K_{1,4}$. It remains to realize the dihedral group $D_4$ of order 8. 
Let $K_1\vee(K_2\sqcup2K_1)$ the graph obtained from $K_2\sqcup2K_1$ by adjoining one new vertex adjacent to every existing vertex (Figure \ref{fig:graph-H}). 
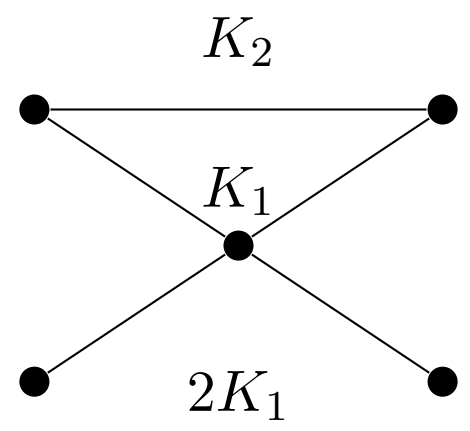
\begin{figure}[ht]
\centering
\begin{tikzpicture}[
    scale=0.9,
    every node/.style={circle,fill=black,inner sep=2pt}
]
% Universal vertex
\node (u) at (0,0) {};

% Vertices of K_2
\node (a) at (-1.5,1) {};
\node (b) at (1.5,1) {};

% Two isolated vertices before taking the join
\node (c) at (-1.5,-1) {};
\node (d) at (1.5,-1) {};

% Edge of K_2
\draw (a)--(b);

% Edges from the universal vertex
\draw (u)--(a);
\draw (u)--(b);
\draw (u)--(c);
\draw (u)--(d);

\node[fill=none] at (0,1.5) {$K_2$};
\node[fill=none] at (0,0.4) {$K_1$};
\node[fill=none] at (0,-1.1) {$2K_1$};

\end{tikzpicture}
\caption{The graph $K_1\vee(K_2\sqcup2K_1)$.}
\label{fig:graph-H}
\end{figure}

A direct computation gives
\[
D(K_1\vee(K_2\sqcup2K_1),x)=x(x^4+5x^3+8x^2+4x+1).
\]
Using MAGMA we can verify that the quartic factor is irreducible with Galois group $D_4$. Hence every transitive subgroup of $S_4$ except $A_4$ is realized.
\end{proof}

Thus $A_4$ is the first even-order transitive group of degree at most four for which we do not know a realization by a domination polynomial.
In degree five, two transitive groups are presently realized in our examples: the cycle $C_8$ has
\[
 D(C_8,x)=x^3(x^5+8x^4+28x^3+48x^2+38x+8),
\]
with quintic Galois group $S_5$ and the graph $F$ used in Remark \ref{rmk:D5} has an irreducible quintic factor with Galois group $D_5$. We know of no domination realizations of the remaining transitive quintic groups $\Z/5\Z$, $F_{20}$, and $A_5$; the first would be excluded by Oboudi's conjecture.

The two viewpoints above give complementary information. Table \ref{tab:small-groups} summarizes the known realizations and open cases for the chromatic and domination inverse Galois problems among groups of order at most eight. 

\begin{table}[htbp]
\centering
\begin{tabular}{|ccc|}
\hline 
Group & Chromatic & Domination\\
\hline 
$\Z/2\Z$ & yes ($C_4$) & yes ($K_3$, $\dagger$)  \\
$\Z/3\Z$ & yes (Morgan \cite[Theorem 4]{Morgan2012}) & conditionally no\\
$(\Z/2\Z)^2$ & yes ($C_9$) & yes ($K_8$, $\dagger$)\\
$\Z/4\Z$ & yes ($C_6$) & yes ($K_5$, $\dagger$)\\
$\Z/5\Z$ & yes (Gilmore \cite[Proposition 32]{DelbourgoMorgan2019}) & conditionally no\\
$S_3$ & yes ($\theta_{2,2,2}$, \cite[Table B.1]{Morgan2012}) & yes ($K_{1,3}$, Theorem \ref{thm:stars})\\
$\Z/6\Z$ & yes ($C_8$) & yes ($K_7$, $\dagger$)\\
$\Z/7\Z$ & --- &  conditionally no\\
$\Z/4\Z \times \Z/2\Z$ & yes ($C_{17}$) & yes ($K_{16}$ , $\dagger$)\\
$(\Z/2\Z)^3$ & yes ($C_{25}$) & yes ($K_{24}$ , $\dagger$)\\
$D_4$ & yes (Morgan \cite[Theorem 7]{Morgan2012}) & yes ($K_1\vee(K_2\sqcup2K_1)$, Figure \ref{fig:graph-H}) \\
$Q_8$ & --- & --- \\
$\Z/8\Z$ & --- & --- \\
\hline 
\end{tabular}
\caption{\small A dagger ($\dagger$) indicates a realization obtained via Proposition~\ref{prop:complete}; a dash (---) indicates that no realization is currently known, while “conditionally no” indicates that nonexistence follows from Oboudi's conjecture via Corollary~\ref{cor:odd-oboudi}.}
\label{tab:small-groups}
\end{table}

Thus, the only unresolved cases of even order at most eight are the cyclic group $\Z/8\Z$ and the quaternion group $Q_8$. The same two groups are also unresolved in the chromatic comparison through order eight. The odd cyclic case $\Z/7\Z$ is included in Morgan's question on $\Z/n\Z$ for odd $n>5$~\cite[Problem~2]{Morgan2012}.

\begin{question}\label{q:order-eight}
Do $\Z/8\Z$ and $Q_8$ occur as Galois groups of chromatic polynomials of connected graphs? Do they occur as Galois groups of domination polynomials of connected graphs?
\end{question}

\subsection*{Acknowledgments}
The authors thank Mario Pineda Ruelas for his insightful questions.
The second author was partially supported by the SECIHTI Grants: CBF 2023-2024-224 and CBF-2025-I-797.

\subsection*{Declaration on the use of AI} 
The authors used ChatGPT 5.6 Sol (OpenAI) for language editing, bibliographic assistance, and exploratory computations concerning Galois groups of domination polynomials. The authors verified the mathematical content and references, and the computational results were verified using MAGMA. The authors take full responsibility for the content of the paper.

\textsc{.\\ Universidad Autonoma Metropolitana - Unidad Iztapalapa \\ Departamento de Matemáticas, Edificio Anexo al T \\
Av San Rafael Atlixco No.186 Col.Vicentina \\
C.P. 09340, Iztapalapa, Mexico City, Mexico.}\\
E-mail: \texttt{coal@xanum.uam.mx}\\
E-mail: \texttt{adrian.zenteno@xanum.uam.mx}

\end{document}